\documentclass[preprint,12pt]{elsarticle}

\usepackage{amssymb}
\usepackage[cp1252]{inputenc}
\usepackage{graphicx}
\usepackage{theorem}

\theorembodyfont{\normalfont}

\newtheorem{defi}{Definition}[section]
\newtheorem{teo}[defi]{Theorem}
\newtheorem{lem}[defi]{Lemma}
\newtheorem{prop}[defi]{Proposition}
\newtheorem{cor}[defi]{Corollary}
\newtheorem{obs}[defi]{Remark}
\newtheorem{ejem}[defi]{Example}

\newcommand{\N}{\mathbb{N}}
\newcommand{\R}{\mathbb{R}}
\newcommand{\K}{\mathbb{K}}
\newcommand{\C}{\mathbb{C}}

\newcommand{\Bo}{\mathbb{B}}
\newcommand{\Boloc}{\Bo(\mu, X)^{loc}}
\newcommand{\Pet}{\mathbb{P}}
\newcommand{\Du}{\mathbb{D}}
\newcommand{\Petloc}{\Pet(\mu, X)^{loc}}

\newcommand{\Ri}{\mathcal{R}}

\newcommand{\No}{\mathcal{N}_0}

\newcommand{\Lw}{L^1_w(\nu)}
\newcommand{\Lf}{L^1(\nu)}

\newcommand{\AR}{\mathcal{R}}
\newcommand{\ARl}{\mathcal{R}^{loc}}
\newcommand{\M}{L^0(\mathcal{R}^{loc})}
\newcommand{\Bx}{B_{X^\ast}}
\newcommand{\xe}{x^\ast}
\newcommand{\Sigf}{\Sigma^f}
\newcommand{\ds}{\displaystyle}
\newcommand{\re}{\emph}
\newcommand{\proof}{\textbf{Proof. }}

\journal{Indagationes Mathematicae}

\begin{document}

\begin{frontmatter}



\title{$L^1$-spaces of  vector measures with vector density}


\author{Celia Avalos-Ramos}

\address{C.U.C.E.I. Universidad de Guadalajara\\
Blvd. Gral. Marcelino Garc\'ia Barrag\'an 1421, Col.  Ol\'impica\\ 
C.P.  44430, Guadalajara, Jal., M\'exico.\\
e-mail: celia.avalos@academicos.udg.mx}

\begin{abstract}
Let $F$ be a function with values in a Banach space. When $F$ is locally (Pettis or Bochner) integrable with respect to a locally determined positive measure, 
a vector measure $\nu_F$ with density $F$ defined on a $\delta$-ring is obtained.   
 We present the existing connection between the spaces $L^1_w(\nu_F)$, $L^1(\nu_F)$ and $L^1(|\nu_F|)$ and the spaces of Dunford, Pettis or Bochner integrable functions.
\end{abstract}


\end{frontmatter}



\section{Introduction}

Let us consider a Banach space $X$ and a $\sigma$-algebra $\Sigma$ on $\Omega$. Since the integrability of $X$-valued functions with respect to a positive finite measure defined in $\Sigma$  was introduced by B. J. Pettis and S. Bochner in the thirties of the last century  this theory has been well studied by several authors; see \cite{Diestel}, \cite{Lang}, \cite{Stefansson}.
After that  the theory of integration of scalar functions  with respect to  $X$-valued measures defined on $\Sigma$ was created by R. G. Bartle, N. Dunford and J. Schwartz in \cite{Bartle}.
It turns out that the Pettis and the Bochner integrals define vector measures, these ones were studied among other by J. Diestel and J. J. Uhl  in \cite{Diestel}. 
On the other hand  N. Dinculeanu and J. J. Uhl  
introduced the concept of $\Ri$-locally (Pettis and Bochner) integrable function with respect to a locally $\sigma$-finite positive measure $\lambda$ defined on a $\delta$-ring $\AR$; see  \cite{DinculeanuUhl}. 
The extension of the integration theory with respect to vector measures defined on $\delta$-rings was developed by D. R. Lewis \cite{Lewis},   P. R. Masani and H. Niemi  \cite{Masani1} and O. Delgado \cite{Delgado}. As one might expect the $\AR$-locally (Pettis and Bochner) integrals define vector measures on $\AR$.


 In this paper we consider a locally determined positive measure $\mu$ on $\Sigma$. Then  
  the collection $\Sigf$ consisting of those subsets in $\Sigma$ which have $\mu$-finite measure is a $\delta$-ring and if $\mu$ is restricted to $\Sigf$, a locally $\sigma$-finite measure is obtained; so makes sense to consider vector functions which are $\Sigf$-locally (Pettis or Bochner) integrable, these functions will be simply called \re{locally (Pettis  or Bochner) integrable functions}. In section \ref{BochPett}  we recall the basic concepts relative to the Dunford, the Pettis and the Bochner integrals with respect to a positive measure and also the main results about vector measures defined on $\delta$-rings. In section \ref{vmvd} we study briefly the vector measure $\nu_F$ defined on $\Sigf$ by the integral of a locally Pettis integrable function $F:\Omega\rightarrow X$ over each $B\in\Sigf$ and a description of its corresponding semivariation is given. If additionally  the function $F$ is locally Bochner integrable we provide a characterization of the variation of the measure $\nu_F$. Finally in section \ref{snuF} we present the existing connection between the weakly $\nu_F$-integrable, $\nu_F$-integrable  or $|\nu_F|$-integrable functions and the Dunford, Pettis or Bochner integrable functions as well as their corresponding integrals. In this way Theorems $8$ and $13$ established by G.F. Stefansson in \cite{Stefansson} are generalized in two directions, namely the positive measure is no longer necessarily finite but locally determined and the vector function $F$ is now locally integrable. Besides we obtain some conditions to determine when a locally Pettis integrable function is in fact Pettis integrable.

\section{Preliminary results}\label{BochPett}

\subsection{Bochner and Pettis integrals} 

Throughout the paper $\Omega$ will be a non empty set and  $X$ stands for a Banach space over $\K$ ($\R$ or $\C$). We denote by $X^*$ and $B_X$  its dual space and  its unit ball respectively.  Let us consider a $\sigma$-algebra  $\Sigma$ on $\Omega$ and a positive measure $\mu:\Sigma\rightarrow [0,\infty]$. The set $\Sigf$ consists of the subsets  $B\in\Sigma$ such that $\mu(B)<\infty$ and $\No(\mu)$ is the collection of $\mu$-null sets. Recall that the measure $\mu$ is said to be \re{semi-finite} if for each set $A\in\Sigma$ such that $\mu(A)>0$, there exists a subset $B\in\Sigf$ satisfying that $B\subset A$ and $0<\mu(B)$. The measure $\mu$ is \re{locally determined} if it is semi-finite and $\Sigma=\{A\subset\Omega \ | \ A\cap B\in\Sigma, \ \forall \ B\in\Sigf\}$.

We denote by $St(\mu,X)$ the vector space of $X$-valued simple functions whose support has finite measure. An $X$-valued function $F:\Omega\rightarrow X$ is said to be \re{strongly $\mu$-measurable} if there exists a sequence  $\{S_n\}\subset St(\mu,X)$, which converges pointwise to $F$ $\mu$-a.e. and to be \re{weakly $\mu$-measurable} if $\langle F,x^*\rangle:\Omega\rightarrow\R$ is strongly $\mu$-measurable for any $x^*\in X^*$. Clearly each strongly $\mu$-measurable function is a weakly $\mu$-measurable function.  We say that two functions $F,G:\Omega\rightarrow X$ are \re{weakly equal} $\mu$-a.e. if $\langle F,x^*\rangle=\langle G,x^*\rangle$ $\mu$-a.e. for all $x^*\in X^*$. 
 We will denote by $L^0(\mu,X)$  the vector space that consists of the equivalence classes that are obtained by identifying strongly $\mu$-measurable functions if they are equal $\mu$- a.e.  and $L^0_w(\mu,X)$ the vector spaces formed by the equivalence classes that are obtained when we identified weakly $\mu$-measurable functions if they are weakly equal $\mu$-a.e. We write $\Bo(\mu,X)$ to indicate the Banach space of the \re{Bochner integrable} functions, namely the functions $F\in L^0(\mu,X)$ such that $\|F\|_X\in L^1(\mu)$, with the norm defined by $\|F\|_1=\int_\Omega \|F\|_Xd\mu$. On the other hand a function $F\in L^0_w(\mu,X)$ is \re{Dunford integrable} when $\langle F,x^*\rangle\in L^1(\mu)$, $\forall \ x^*\in X^*$, if additionally for each $A\in\Sigma$ there exists a vector $x_A\in X$ such that $\int_A \langle F,x^*\rangle d\mu=\langle x_A,x^*\rangle$; $\forall \ x^*\in X^*$, the function $F$ is \re{Pettis integrable} and  the vector $x_A$ is called  the \re{Pettis integral} of $F$ over $A$ and it is denoted by $\Pet-\int_A F d\mu$. We write $\Du(\mu,X)$ and $\Pet(\mu,X)$ for the vector spaces consisting of the Dunford and Pettis integrable functions respectively.
 \begin{lem}
   The space $\Du(\mu,X)$ is a normed  space with the norm given by 
   \[
    \|F\|_\Du=\sup_{x^*\in B_{X^*}}\int_\Omega |\langle F,x^*\rangle|d\mu, \ \forall \ F\in \Du(\mu,X).
  \] 
 \end{lem}
 \proof
    Observe that to obtain the conclusion it only remains to establish that $\|F\|_\Du<\infty$, $\forall \ F\in\Du(\mu,X)$. So, let us fix $F\in\Du(\mu,X)$ and define  $T:X^\ast\rightarrow L^1(\mu)$ by   $T(x^\ast)=\langle F,x^\ast\rangle, \ \forall \ x^\ast\in X^\ast$.
 Clearly $T$ is a well defined linear operator. Now take $\{\xe_n\}\subset X^\ast$ and $\xe\in X^\ast$ such that $\xe_n\rightarrow x$. Let us assume that there exists $g\in L^1(\mu)$ satisfying that $T\xe_n\rightarrow g$ in $L^1(\mu)$. Proceeding as in \cite[p.46]{DelgadoT} we get a  subsequence  $\{\xe_{n_k}\}\subset\{\xe_n\}$ such that  $\langle F,\xe_{n_k}\rangle= T\xe_{n_k}\rightarrow g$ $\mu$-a.e.  On the other hand $\langle F(t),\xe_n\rangle\rightarrow \langle F(t),\xe\rangle$, $\forall \ t\in\Omega$. Thus, $T(\xe)=\langle F(t),\xe\rangle=g$, $\mu$-a.e. By the closed graph theorem we get that $T$ is bounded. Therefore 
\[
 \|F\|_\Du=\sup_{x^\ast\in B_{X^\ast}}\int_\Omega |\langle F,x^\ast\rangle|d\mu=\sup_{x^\ast\in B_{X^\ast}}\|T\xe\|_{L^1(\mu)}\leq \|T\|. \ \blacksquare
\]
 Since $\Pet(\mu,X)\subset\Du(\mu,X)$ it turns out that $\Pet(\mu,X)$ is also a normed space with the same norm $\|\cdot\|_\Du$ which will be denoted by $\|\cdot\|_\Pet$ in this case.   It is well know that $\Bo(\mu,X)\subset\Pet(\mu,X)$ with $\|F\|_\Pet\leq\|F\|_1$ and $\Bo-\int_A F d\mu=\Pet-\int_A Fd\mu$, $\forall \ F\in\Bo(\mu,X)$.

\subsection{Integration with respect to measures defined on $\delta$-rings}

A family $\AR$ of subsets of $\Omega$ is a \re{$\delta$-ring} if $\AR$ is a ring which is closed under countable intersections. 
We denote by $\ARl$ the $\sigma$-algebra of all sets $A\subset\Omega$ such that $A\cap B\in\AR$, $\forall \ B\in\AR$. Given $A\in\ARl$ we indicate by $\AR_A$ the $\delta$-ring $\{B\subset A:B\in\AR\}$ and by $\pi_A$  the collection of finite families of pairwise disjoint sets in $\AR_A$. Note that if $\Omega\in\AR$, then $\AR$ is a $\sigma$-algebra, and in this case we have that $\ARl=\AR$. Moreover, for each $B\in\AR$ it turns out that $\AR_B$ is a $\sigma$-algebra.\vspace{5 pt}

 A \re{scalar measure} is a function $\lambda:\AR\rightarrow \K$  satisfying that if $\{B_n\}\subset \AR$, is a family of pairwise disjoint sets such that $\bigcup_{n=1}^{\infty}B_n\in \AR$, then $\sum_{n=1}^{\infty} \lambda(B_n)=\lambda\left(\bigcup_{n=1}^{\infty}B_n\right)$. The \re{variation of $\lambda$} is the countably additive measure $|\lambda|:\ARl\rightarrow[0,\infty]$ defined by
\[
   |\lambda|(A):=\sup\left\{\sum_{j=1}^{n}|\lambda(A_j)|
     :\{A_j\}\in\pi_A\right\}.
\]
A function $f\in\M$ is \re{$\lambda$-integrable} if $f\in L^1(|\lambda|)$. We denote by $L^1(\lambda)$ the vector space consisting of the equivalence classes of $\lambda$-integrable functions when we identify two functions if they are equal $|\lambda|$-a.e.  
\vspace{5 pt}

 A set function $\nu:\AR\rightarrow X$ is a \re{vector measure} if  for any collection $\{B_n\}\subset\AR$  of pairwise disjoint sets satisfying that $\bigcup_{n=1}^\infty B_n \in \AR$ we have  $\sum_{n=1}^\infty \nu(B_n)=\nu(\bigcup_{n=1}^\infty B_n)$. A vector measure $\nu$ is called \re{strongly additive} if $\nu(B_n)\rightarrow 0$ whenever $\{B_n\}$ is a disjoint sequence in $\AR$.  The \re{variation of $\nu$} is the positive measure $|\nu|$ defined in $\ARl$ by 
 \[
 |\nu|(A):=\sup\left\{\sum_{j=1}^n \|\nu(A_j)\|_X : \{A_j\}\in\pi_A \right\}.
 \]
The \re{semivariation of $\nu$} is the function $\|\nu\|:\ARl\rightarrow [0,\infty]$ given by
\[
  \ds{\|\nu\|(A):=\sup\{|\langle\nu,x^*\rangle|(A):x^*\in B_{X^*}\}},
\]
where $|\langle\nu,x^*\rangle|$ is the variation of the scalar measure  $\langle\nu,x^*\rangle:\AR\rightarrow \K$, defined by
\[
  \langle\nu,\xe\rangle(B)=\langle\nu(B),\xe\rangle, \ \forall \ B\in\AR.
\]
The semivariation of $\nu$ is finite in $\AR$ and for any $A\in\ARl$ satisfies $\|\nu\|(A)\leq |\nu|(A)$. A set $A\in\ARl$ is said to be $\nu$-null if $\|\nu\|(A)=0$. We will denote by $\No(\nu)$ the collection of $\nu$-null sets. It turns out that $\No(\nu)=\No(|\nu|)$. Moreover $A\in\No(\nu)$ if and only if $\nu(B)=0$, $\forall \ B\in\AR_A$. We say that two functions $f,g\in\M$ are equal $\nu$-a.e. if they are equal outside of a set in $\No(\nu)$.
We define $L^0(\nu)$ as  the space of equivalence classes of functions in $L^0(\ARl)$, where two functions are identified when they are equal  $\nu$-a.e. 

A function $f\in \M$ is \re{weakly $\nu$-integrable}, if $f\in L^1(\langle\nu,x^*\rangle)$, for each  $x^*\in X^*$. We will denote by $\Lw$ the subspace of $L^0(\nu)$ of all weakly $\nu$-integrable functions. With the norm given by 
\[
\|f\|_\nu:=\sup\left\{\int_\Omega|f|\hspace{3 pt} d|\langle\nu,x^*\rangle|:x^*\in B_{X^*}\right\},
\]
 $\Lw$ is a Banach space. 

A function $f\in\Lw$ is \re{$\nu$-integrable}, if for each $A\in\AR^{loc}$ there exists a vector \mbox{$x_A\in X$}, such that
  \begin{equation}
    \langle x_A,\xe\rangle=\int_A f d\langle\nu,\xe\rangle, \ \forall \ \xe\in X^*.
  \end{equation}
In this case the vector $x_A$ is denoted by $\int_A f d\nu$.  With the norm $\|\cdot\|_\nu$ the subset of all $\nu$-integrable functions is a closed subspace of $\Lw$ and it will be denoted by $\Lf$. Therefore $\Lf$ is also a Banach space . We indicate by $S(\AR)$ the collection of simple functions in $\M$ which have support in $\AR$. It turns out that $S(\AR)$ is a dense subspace of $\Lf$.  Finally the integral operator $I_\nu:\Lf\rightarrow X$ defined by $I_\nu(f)=\int_\Omega f d\nu$, is linear and bounded.

\section{Vector measures with vector density}\label{vmvd}

Recall that $(\Omega,\Sigma,\mu)$ is a positive measure space and $X$ a Banach space.   Given a vector function $F\in\Pet(\mu,X)$ let us define the set function $\tilde{\nu}_F:\Sigma\rightarrow X$ by
\begin{equation}\label{dvmB}
   \tilde{\nu}_F(A)=\Pet-\int_A F d\mu,  \ \ \forall \ A\in\Sigma.
\end{equation}


In the case when $\mu$ is finite it is well known that $\tilde{\nu}_F$ is a vector measure \cite[Thm. II. 3.5]{Diestel}.  The next result generalizes this fact, it can be established in the similar way, using the Orlicz-Pettis Theorem \cite[Cor. I.4.4]{Diestel}.

\begin{teo}\label{nutildeFmv}
The set function $\tilde{\nu}_F$ defined on (\ref{dvmB}) is a vector measure with semivariation 
\begin{equation}\label{semivariacionnutilde}
  \|\tilde{\nu}_F\|(A)=\sup_{x^*\in B_{X^*}}\int_A |\langle F,x^*\rangle|d\mu, \ \forall \ A\in\Sigma.
\end{equation}
\end{teo}
\proof
 Let us fix $x^*\in X^*$ and take a pairwise disjoint countable collection $\{A_n\}\subset\Sigma$, then
 \[
 \begin{array}{lll}
\ds{\left\langle \Pet-\int_{\bigcup_n A_n} F d\mu,x^*\right\rangle}&=&\ds{\int_{\bigcup_n A_n} \langle F,x^*\rangle d\mu} \vspace{3 pt}\\ 
&=&\ds{\sum_{n=1}^\infty \int_{A_n} \langle F,x^*\rangle d\mu=\sum_{n=1}^\infty \left\langle \Pet-\int_{A_n}Fd\mu,x^*\right\rangle}.
 \end{array}
 \]
So, $\tilde{\nu}_F$ is weakly $\sigma$-aditive, by the Orlicz-Pettis Theorem \cite[Cor. I.4.4]{Diestel} $\tilde{\nu}_F$ is a vector measure.
On the other hand since $\langle \tilde{\nu}_F, x^*\rangle(A)=\int_A \langle F,x ^*\rangle$, $\forall$ $A\in \Sigma$, $\mu$ is a positive measure  and $\langle F,x^*\rangle\in L^1(\mu)$ from \cite[Thm. 6.13]{Rudin} we get that
\[
  |\langle \tilde{\nu}_F,x^*\rangle|(A)=\int_A |\langle F,x^*\rangle|d\mu.
\]
It is follows (\ref{semivariacionnutilde}). $\blacksquare$ \vspace{10 pt}

Since $\Bo(\mu,X)\subset\Pet(\mu,X)$ we have the following consequence. In order to get (\ref{variacionnutilde}) we can proceed as in \cite[Thm. II.2.4 iv)]{Diestel}.

\begin{cor}
  Let $F\in \Bo(\mu,X)$. Then $\tilde{\nu}_F$ defined on (\ref{dvmB}) is a vector measure with bounded variation such that
  \begin{equation}\label{variacionnutilde}
      |\tilde{\nu}_F|(A)=\int_A \|F\|_X d\mu, \ \forall \ A\in\Sigma.
  \end{equation}
\end{cor}

\vspace{15 pt}
Hereafter let us consider  a  locally determined positive measure $\mu:\Sigma\rightarrow [0,\infty]$. It turns out that $\Sigf$ is a $\delta$-ring such that $(\Sigf)^{loc}=\Sigma$ and the restriction of $\mu$ to $\Sigf$, denoted by $\lambda$, is a scalar measure such that $|\lambda|=\mu$ (\cite[Lemma 4.3]{AvalosG}).


Now we will study a kind of vector functions that  include vector measures  having Pettis or Bochner functions as density functions. 

\begin{defi} Let $F:\Omega\rightarrow X$ be a vector function.
\begin{enumerate}
\item[a)] The function $F$ is \re{locally Pettis integrable} if $F$ is weakly $\mu$-measurable function and $\chi_B F\in \Pet(\mu,X)$, $\forall$ $B\in\Sigma^f$.   The collection of equivalence classes obtained by identifying  locally Pettis integrable functions if they are weakly equal  $\mu$-a.e.  will be denoted by $\Petloc$.

\item[b)]  $F$ es \re{locally Bochner integrable} if $F$ is strongly $\mu$-measurable function and $\chi_B F\in \Bo(\mu,X)$, $\forall$ $B\in\Sigma^f$.  The collection of equivalence classes obtained by identifying  locally Bochner integrable functions if they are equal  $\mu$-a.e.  will be denoted by $\Boloc$.
\end{enumerate}
\end{defi}

\begin{obs}
Observe that $\Petloc$ and $\Boloc$ are vector spaces. Additionally  we have that 
\[
\Bo(\mu,X)\subset\Boloc\subset\Petloc.
\]
\end{obs}

The following example shows that the containment $\Bo(\mu,X)\subset\Boloc$ can be proper.

\begin{ejem}\label{elejmplo}
Let us fix $x\in X$ and assume that $f:\Omega\rightarrow \R$ is such that  $\chi_B f\in L^1(\mu)$, $\forall \ B\in\Sigf$ (c.f. \cite[Def. 2.14 c)]{Masani1}). Now define $F:\Omega\rightarrow X$ by
\begin{equation}\label{defF}
  F(t):=f(t) x.
\end{equation}

Let us see that $\chi_BF\in\Bo(\mu,X)$, $\forall$ $B\in \Sigma$. Since $f\chi_B\in L^0(\Sigma)$, for all $B\in\Sigf$, it turns out that $f\in L^0(\Sigma)$. Take $\{s_n\}\subset S(\Sigma)$ such that $s_n\rightarrow f$ and $|s_n|\leq |f|$, $\forall$ $n\in\N$. Fix $B\in\Sigf$. For each $n\in\N$, define $S_n:\Omega\rightarrow X$ by $S_n(t)= \chi_B s_n(t)x$. Since $\chi_Bf\in L^1(\mu)$, we have that $\chi_Bs_n\in L^1(\mu)$. And so $S_n\in St(\mu,X)$ and $S_n(t)\rightarrow \chi_B f(t) x=\chi_BF(t)$, $\forall$ $t\in\Omega$ indicating that $\chi_B F$ is strongly $\mu$-measurable.

Further
\[
 \ds{\int_B \|F\|_Xd\mu=\int_B |f|\|x\|_Xd\mu=\left(\int_B |f|d\mu\right) \|x\|_X }.
\]
Since $\chi_B f\in L^1(\mu)$, it follows that  $\ds\int_B \|F\|_Xd\mu<\infty$, and so $\chi_BF\in\Bo(\mu,X)$. Thus $F\in \Boloc$. Finally observe that $F\in \Bo(\mu,X)$ if and only if $f\in L^1(\mu)$.

\end{ejem}

\vspace{10 pt}

As a consequence of Theorem \ref{nutildeFmv} we obtain the following result.

\begin{prop}
Let $F\in\Petloc$. Then the set function $\nu_F:\Sigf\rightarrow X$ defined by
\begin{equation}\label{dmvBr}
\nu_F(B):=\Pet-\int_B Fd \mu, \ \forall \ B\in\Sigf,
\end{equation}
is a vector measure such that
\begin{equation}\label{semivarnuF}
\|\nu_F\|(A)=\sup_{x^*\in B_{X^*}}\int_A |\langle F,x^*\rangle|d\mu, \ \forall \ A\in\Sigma.
\end{equation}
\end{prop}
\proof
Let us show that $\nu_F$ is a vector measure.  Let $\{B_n\}\subset\Sigf$ be a disjoint collection such than $B:=\bigcup_{n=1}^\infty B_n\in\Sigf$. Since $\chi_B F\in \Pet(\mu,X)$, by Theorem \ref{nutildeFmv} we have
\[
\begin{array}{lll}
\nu_F(B)&=& \ds{\Pet-\int_B Fd\mu=\Pet-\int_B \chi_BFd\mu} \vspace{5 pt} \\
                  &=&\ds{\sum_{n=1}^\infty \Pet-\int_{B_n}\chi_B Fd\mu=\sum_{n=1}^\infty \Pet-\int_{B_n} Fd\mu=\nu_F(B_n)}.
\end{array}
\]
Thus $\nu_F$ is a vector measure. Now fix $x^*\in X^*$. Since $F\in\Petloc$ we have $\langle F,x^*\rangle\in L^1_{loc}(\lambda)$. From \cite[Thm. 2.31]{Masani1} we have that the variation of the scalar measure $\mu_{\langle F,x^*\rangle}:\Sigf\rightarrow \K$ defined by $\mu_{\langle F,x^*\rangle}(B)=\int_B \langle F,x^*\rangle d\mu$, $\forall \ B\in\Sigf$ is given by
\[
|\mu_{\langle F,x^*\rangle}|(A)=\int_A |\langle F,x^*\rangle|d\mu, \ \forall \ A\in\Sigma.
\]

On the other hand notice that for each $B\in\Sigf$, we have
\begin{equation}\label{medescnuF}
\ds{\langle \nu_F,x^*\rangle(B)=\left\langle\Pet-\int_B Fd\mu, x^*\right\rangle=\int_B \langle F,x^*\rangle d\mu=\mu_{\langle F,x^*\rangle}(B).}
\end{equation}
Therefore
\begin{equation}\label{varmedescnuF}
|\langle \nu_F,x^*\rangle|(A)=\int_A |\langle F,x^*\rangle|d\mu.
\end{equation}
From this we get (\ref{semivarnuF}). $\blacksquare$

\vspace{10 pt}

Observe that from (\ref{semivarnuF}) we have that if $A\in\No(\mu)$, then $\|\nu_F\|(B)=0$, $\forall$ $B\in\Sigf_A$. So $\No(\mu)\subset\No(\nu_F)$. Also observe that $\|\nu_F\|(B)=\|\chi_BF\|_\Pet$, $\forall $ $B\in\Sigf$. Moreover from the Dunford integrability definition we get our next result.
\begin{cor}
  Let $F\in\Petloc$. Then $F$ is Dunford integrable if and only if $\nu_F$ has bounded semivariation. In particular if $F\in\Pet(\mu,X)$, then $\|\nu_F\|(A)=\|\chi_A F\|_\Pet$.
\end{cor}

In the case that $F$ is locally Bochner integrable we have a useful characterization of the variation of $\nu_F$, as follows.

\begin{prop}\label{varBochner}
If $F\in\Boloc$, then $\ds{|\nu_F|(A)=\int_A \|F\|_X d\mu}$, $\forall$  $A\in\Sigma$.
\end{prop}

\proof
Take $B\in\Sigf$. Notice that $\nu_F(A\cap B)= \nu_{\chi_BF}(A)$, $\forall$ $A\in\Sigma$. Since $\chi_BF\in\Bo(\mu,X)$ we have that
\[
  |\nu_F|(B)=|\nu_{\chi_BF}|(B)=\int_B\|F\|_Xd\mu.
\]
Hence
\begin{equation}\label{varnuF}
|\nu_F|(A)=\sup_{B\in\Sigf_A}|\nu_F|(B)=\sup_{B\in\Sigf_A}\int_B\|F\|_Xd\mu=\int_A\|F\|_Xd\mu, \ \forall \ A\in\Sigma.  \  \  \blacksquare
\end{equation}

\begin{obs}
Let us note that if $F:\Omega\rightarrow X$ is a strong $\mu$-measurable function, by the previous result we obtain that:
\begin{center}
$F\in\Bo(\mu,X)$ if and only if $F\in\Boloc$ and $\nu_F$ has bounded variation.
\end{center}
\end{obs}

\begin{ejem}\label{elejmplo2}
Let us return to the example \ref{elejmplo}. It was shown there that $F$ defined in (\ref{defF}) is locally Bochner integrable.  In particular $F$ is locally Pettis integrable. Let us obtain now the vector measure $\nu_F$,  its variation and its semivariation. See that for each  $B\in\Sigf$

\[
   \nu_F(B)=\Bo-\int_B Fd\mu=\left(\int_B f d\mu\right) x, \ \forall \ B\in\Sigf.
 \]

Now take $A\in\Sigma$, from (\ref{semivarnuF}) and (\ref{varnuF}) 
 \[
 \begin{array}{lll}
    \|\nu_F\|(A)&=&\ds{\sup_{\xe\in\Bx}\int_A|\langle F,\xe\rangle|d\mu
                               =\sup_{\xe\in\Bx}\int_A |f||\langle x,\xe\rangle|d\mu} \vspace{8 pt}\\ 
    							&=&\ds{\left(\int_A |f|d\mu\right)\sup_{\xe\in\Bx}|\langle x,\xe\rangle|
    							=\left(\int_A |f|d\mu\right)\|x\|_X}\vspace{8 pt}\\
    							&=&\ds{\int_A \|F\|_Xd\mu=|\nu_F|(A)}.
 \end{array}
 \]
Therefore $\|\nu_F\|=|\nu_F|$ in this case.  \vspace{12 pt}
\end{ejem}

\section{The space of $\nu_F$-integrable functions}\label{snuF}

 When $F$ is a locally Pettis or Bochner integrable function we have constructed  the vector measure $\nu_F$ defined on the $\delta$-ring $\Sigf$. In the present section we will study the spaces $L^0(\nu_F)$, $L^1_w(\nu_F)$, $L^1(\nu_F)$ and $L^1(|\nu_F)|)$ associated to this vector measure through the operator $M_F$ which to each function $g$ assigns the function $gF$.  The following lemmas allow us to conclude that $M_F:L^0(\nu_F)\rightarrow L^0_w(\mu,X)$ or $L^0(\mu,X)$ is well defined. Clearly $M_F$ is a linear operator.

\begin{lem}\label{gFescmumed}
Let $F:\Omega\rightarrow X$ be a function and $g\in L^0(\Sigma)$. 
\begin{enumerate}
 \item[i)]  If $F$ is strongly $\mu$-measurable, then  $gF$ is strongly $\mu$-measurable. 
 \item[ii)] If $F$ is weakly $\mu$-measurable, then $gF$ is weakly $\mu$-measurable.
\end{enumerate}   
\end{lem}
\proof
i) Observe that if $\varphi\in S(\Sigma)$ and $S\in St(\mu,X)$, then $\varphi S\in St(\mu,X)$. Let us assume that $F$ is strongly $\mu$-measurable. Take $\{\varphi_n\}\subset S(\Sigma)$ and $\{S_n\}\subset St(\mu,\Sigma)$ such that $\varphi_n\rightarrow g$ and $S_n\rightarrow F$, $\mu$-a.e. Thus $\varphi_nS_n\in St(\mu,X)$, $\forall \ n\in\N$ and $\varphi_nS_n\rightarrow gF$, $\mu$-a.e. It follows that $gF$ is strongly $\mu$-measurable.\vspace{5 pt} \\
ii) By definition if $F$ is weakly $\mu$-measurable, we have that for each $x^*\in X^*$,  the function $\langle F,x^*\rangle$ is strongly $\mu$-measurable.  Using i) we obtain that $\langle gF,x^*\rangle=g\langle F,x^*\rangle\in L^0(\mu,X)$ $\forall \ x^*\in X^*$. 
$\blacksquare$

\vspace{10 pt}

\begin{lem}\label{esclamctp}
 Let $F\in\Petloc$, $\{g_n\}\subset L^0(\Sigma)$ and $g,h\in L^0(\Sigma)$.
 \begin{enumerate}
 \item[i)]   If $g=h$, $\nu_F$-a.e., then $gF=hF$, weakly $\mu$-a.e.
 \item[ii)]  If $g_n\rightarrow g$, $\nu_F$-a.e., then $\langle g_nF, x^*\rangle\rightarrow \langle gF,x^*\rangle$, $\mu$-a.e., $\forall \ x^*\in X^*$.
\end{enumerate}    
\end{lem}
\proof
i) Choose $N\in\No(\nu_F)$ such that $g(t)=h(t)$, $\forall \ t\in N^c$. Then $g\chi_{N^c} F=h\chi_{N^c} F$; moreover  $\chi_NF=0$ weakly $\mu$ -a.e. implies that $g\chi_N F=h\chi_N F=0$,  weakly $\mu$-a.e. Thus $gF=hF$,  weakly $\mu$-a.e. \\
ii) Let $N\in\No(\nu_F)$ such that $g_n(t)\rightarrow g(t)$, $\forall \ t\in N^c$. So $g_n\chi_{N^c}\rightarrow g\chi_{N^c}F$ and $g_n\chi_NF=g\chi_NF=0$, weakly $\mu$-a.e. Then for each $x^*\in X^*$, $\langle g_n\chi_{N^c}F,x^*\rangle\rightarrow \langle g\chi_{N^c}F,x^*\rangle$ and $\langle g_n\chi_NF,x^*\rangle=\langle g\chi_NF,x^*\rangle=0$, $\mu$-a.e. Therefore $\langle g_nF,x^*\rangle\rightarrow\langle gF,x^*\rangle$, $\mu$-a.e., $\forall \ x^*\in X^*$. $\blacksquare$

\vspace{15 pt}

Proposition \cite[Prop.8]{Stefansson}  established by G. F. Stefansson for the case  that $F\in \Pet(\mu,X)$ and $\mu$ is a finite positive measure defined on a $\sigma$-algebra is generalized in the next theorem. 

\begin{teo}\label{isoMF}
For $F\in\Petloc$ and $g\in L^0(\Sigma)$, we have that
\begin{enumerate}
  \item[i)]  $g\in L^1_w(\nu_F)$ if, and only if, $gF\in \Du(\mu,X)$. Moreover the restriction to $L^1_w(\nu_F)$ of the operator $M_F$ is a linear isometry from $L^1_w(\nu_F)$ into $\Du(\mu,X)$.
\item[ii)] $g\in L^1(\nu_F)$ if, and only if, $gF\in \Pet(\mu,X)$. Moreover  $M_F:L^1(\nu_F)\rightarrow \Pet(\mu,X)$, the restriction of the operador $M_F$, is a linear isometry such that $I_{\nu_F}=I_\Pet\circ M_F$.
\end{enumerate}    
\end{teo}
\proof
Fix $\xe\in X^\ast$ and consider $\ds{s=\sum_{j=1}^n a_j\chi_{A_j}\in S(\Sigma)}$.  By hypothesis $\chi_BF\in\Pet(\mu,X)$, $\forall \ B\in\Sigf$. It follows that  $sF\in\Pet(\mu,X)$. From (\ref{varmedescnuF}) we obtain
\begin{equation}\label{pasimp}
\begin{array}{lll}
       \ds{\int_\Omega |s| d|\langle\nu_F,\xe\rangle|}&=&\ds{\sum_{j=1}^n} |a_j|\hspace{1 pt} |\langle\nu_F,\xe\rangle|(A_j) 
			=\ds{\sum_{j=1}^n}| a_j|\int_{A_j}|\langle F,\xe\rangle| d\mu \vspace{5 pt} \\ 
			&=&\ds{\int_\Omega\sum_{j=1}^n}| a_j|\chi_{A_j}|\langle F,\xe\rangle| d\mu
			=\int_\Omega |\langle sF,\xe\rangle| d\mu.
\end{array}
\end{equation}
Thus $s\in L^1_w(\nu_F)$ if and only if $sF\in \Du(\mu,X)$. \vspace{5 pt}

Proceeding in the same way, it follows  from  (\ref{medescnuF}) that
 \begin{equation}\label{intmedescnuF}    
     \ds{\int_\Omega s d\langle\nu_F,\xe\rangle}=\int_\Omega \langle sF,\xe\rangle d\mu.
 \end{equation}

Now take $g\in L^0(\nu_F)^+$ and $\{s_n\}\subset S(\Sigma)$ such that $0\leq s_n\uparrow g$, $\nu_F$-a.e. From lemma \ref{esclamctp} we obtain $\langle s_nF,\xe\rangle\rightarrow\langle gF,\xe\rangle$, $\mu$-a.e. Then, $|\langle s_nF,\xe\rangle|\uparrow |\langle gF,\xe\rangle|$, $\mu$-a.e.  By the Monotone Convergence Theorem and (\ref{pasimp}) it turns out that
\begin{equation}\label{palanor}
 \begin{array}{lll}
 |\langle gF,\xe\rangle|d\mu&=&\ds{\lim_{n\rightarrow\infty} \int_\Omega|\langle s_nF,\xe\rangle|d\mu} \vspace{3 pt}\\ 
 &=&\ds{\lim_{n\rightarrow\infty}\int_\Omega s_n d|\langle \nu_F,\xe\rangle|=\int_\Omega|g|d|\langle \nu_F,\xe\rangle|},
 \end{array} 
\end{equation}
showing that $g\in L^1_w(\nu_F)$ if and only if $gF\in \Du(\mu,X)$.   \vspace{8 pt}

By the Dominate Convergence Theorem and (\ref{medescnuF})
\[
   \int_\Omega \langle gF,\xe\rangle d\mu=\lim_{n\rightarrow\infty}\int_\Omega \langle s_nF,\xe\rangle d\mu=\lim_{n\rightarrow\infty}\int_\Omega s_n d\langle \nu_F,\xe\rangle=\int_\Omega g d\langle\nu_F,\xe\rangle.
\]
We conclude from here that $g\in L^1(\nu_F)$ if and only if  $gF\in\Pet(\mu,X)$ and
\begin{equation}\label{opintnuF}
  \int_\Omega gd\nu_F=\Pet-\int_\Omega gFd\mu. 
\end{equation}

Since the involved sets are vector spaces and each $g\in L^0(\nu_F)$ is a linear combination of non negative functions, we obtain the first part in i) and ii).\vspace{5 pt}

Finally take $g\in L^1_w(\nu_F)$, since $|g|\geq 0$ we obtain equality  (\ref{palanor}) with a sequence $\{s_n\}\subset S(\Sigma)$ such that $0\leq s_n\uparrow |g|$, $\nu_F$-a.e. Taking the supremum over $\xe\in\Bx$ it turns out that $\|g\|_{\nu_F}=\|gF\|_\Du$. That is,  $M_F$ restricted to $L^1_w(\nu_F)$ is a linear isometry. Since $L^1(\nu_F)$ and $\Pet(\mu,X)$ are subspaces of $L^1_w(\nu_F)$ and $\Du(\mu,X)$, respectively we conclude that $M_F$ restricted to $L^1(\nu_F)$ is also an isometry. Moreover,  from (\ref{opintnuF}) it follows that $I_{\nu_F}=I_\Pet\circ M_F$.  $\blacksquare$
\vspace{12 pt}

\begin{cor}
Let $F\in L_w^0(\mu,X)$. Then $F\in\Pet(\mu,X)$ if and only if $F\in\Petloc$ and $\nu_F$ is strongly additive. 
\end{cor}
\proof
Let assume that  $F\in\Pet(\mu,X)$. Consider the vector measure $\tilde{\nu}_F:\ARl\rightarrow X$ defined in (\ref{dvmB}). Since $\Sigma$ is a $\sigma$-algebra, it turns out that $\tilde{\nu}_F$ is strongly additive. Observe that $\nu_F$ is the restriction of $\tilde{\nu}_F$ to $\Sigf$, so it follows that it is strongly additive. \vspace{5 pt}

Now assume that $F\in\Petloc$ and that $\nu_F$ is strongly additive. From \cite[Cor. 3.2]{Delgado}  we obtain that $\chi_\Omega\in L^1(\nu_F)$. So, by the previous theorem $F=\chi_\Omega F\in \Pet(\mu,X)$. $\blacksquare$

\vspace{12 pt}

\begin{cor}
  Let $F\in\Petloc$. If $X$ does not contain any subspace isomorphic to $c_0$ and $\nu_F$ is bounded, then $F\in\Pet(\mu,X)$.
\end{cor}
\proof
Since $X$ does not contain any subspace isomorphic to $c_0$ and $\nu_F$ is bounded it turns out that $\nu_F$ is strongly additive \cite[p. 36]{DelgadoT}. Then by the previous corollary $F\in\Pet(\mu,X)$. $\blacksquare$
 \vspace{15 pt}

The following result gives us the connection between the spaces $L^1(|\nu_F|)$ and $\Bo(\mu,X)$ through the operator $M_F$ in case that $F\in \Boloc$. We will show that, as it occur when $F\in\Petloc$, $M_F$ is a linear isometry in this case.

\begin{prop}\label{isoMFBo}
 Consider $F\in\Boloc$ and $g\in L^0(|\nu_F|)$. Then $g\in L^1(|\nu_F|)$ if and only if $gF\in\Bo(\mu,X)$. Moreover $M_F:L^1(|\nu_F|)\rightarrow\Bo(\mu,X)$  is a linear isometry such that  $I_{\nu_F}(g)=I_\Bo\circ M_F(g)$, $\forall \ g\in L^1(|\nu_F|)$.
\end{prop}
\proof
Clearly $M_F$ is a linear operator, we will see that its image is a subset of $\Bo(\mu,X)$.  
By lemma \ref{esclamctp} we have that the restriction $M_F:L^1(|\nu_F|)\rightarrow \Bo(\mu,X)$ is well defined.

Since the norms in $L^1(|\nu_F|)$ and $\Bo(\mu,X)$ are different from those in $L^1(\nu_F)$ and $\Pet(\mu,X)$ respectively, we need to establish that, under these norms, $M_F$ is also an isometry. \vspace{5 pt}

By hypothesis $F\in\Boloc$, then from (\ref{varnuF}) it follows that $|\nu_F|(B)<\infty$, $\forall \ B\in\Sigf$. So, $S(\Sigf)\subset L^1(|\nu_F|)$. 
Further for each $s=\sum_{j=1}^n a_j\chi_{A_j}\in S(\Sigma)$ we have that 
\[
\begin{array}{lll}
  \ds{\int_\Omega |s|d|\nu_F|}&=&\ds{\int_\Omega \sum_{j=1}^n|a_j|\hspace{1 pt}|\nu_F|(A_j)=\sum_{j=1}^n|a_j|\int_{A_j}\|F\|_Xd\mu}\vspace{5 pt} \\
                          &=&\ds{\int_\Omega \sum_{j=1}^n|a_j|\chi_{A_j}\|F\|_Xd\mu=\int_\Omega \|sF\|_Xd\mu}.
\end{array}
\]

Therefore  $s\in L^1(|\nu_F|)$ if and only if $sF\in\Bo(\mu,X)$. Now consider $g\in L^0(|\nu_F|)$ and take $\{s_n\}\subset S(\Sigma)$ such that $0\leq s_n\uparrow |g|$, $\nu_F$-a.e.  Then $\|s_nF\|_X\uparrow\|gF\|_X$, $\mu$-a.e. By the Monotone Convergence Theorem 
\[
   \int_\Omega \|gF\|_Xd|\mu=\lim_{n\rightarrow\infty}\int_\Omega\|s_nF\|_Xd\mu=\lim_{n\rightarrow\infty}\int_\Omega|s_n|d |\nu_F|=\int_\Omega|g|d|\nu_F|.
\]
Thus we have that $gF\in\Bo(\mu,X)$ if and only if $g\in L^1(|\nu_F|)$. Moreover $\|g\|_{|\nu_F|}=\|gF\|_1$. \vspace{8 pt}

The equality between the operators follows from Proposition \ref{isoMF}. $\blacksquare$ \vspace{12 pt}

\begin{ejem}
 Consider again the function $F$ defined in (\ref{defF}). As we see in Example \ref{elejmplo} $F\in \Boloc$.
 
Take $g\in L^1(\nu_F)$,  from Proposition \ref{isoMF}
\[
  \int_\Omega gd\nu_F=\Pet-\int_\Omega gFd\mu=\left(\int_\Omega gf d\mu\right)x,  
\]
then $gf\in L^1(\mu)$. And so,
\begin{equation}\label{norF}
   \int_\Omega \|gF\|_Xd\mu=\int_\Omega |gf|\|x\|_Xd\mu<\infty.
 \end{equation}
By Lemma \ref{gFescmumed}  $gF$ is strongly  $\mu$-measurable. Thus we have that $gF\in\Bo(\mu,X)$ and by Proposition \ref{isoMFBo},  $g\in L^1(|\nu_F|)$. We conclude that $L^1(|\nu_F|)=L^1(\nu_F)$. And from \cite[Prop. 5.4]{Calabuig2} it follows that $L^1(|\nu_F|)=L^1(\nu_F)=L^1_w(\nu_F)$ . \vspace{12 pt}

\end{ejem}

\subsection*{Acknowledgment}
Many thanks to  Professor F. Galaz-Fontes for the careful reading and valuable comments.





\end{document}